%% file: project1.tex
\documentclass[11pt]{article}

\usepackage[english]{style}
\usepackage{latexsym,amsmath,amssymb,amsfonts,graphics,amscd}
\usepackage[all]{xy}
\usepackage{pgf,tikz}
\usetikzlibrary{arrows}
\usepackage{amssymb}
\input{defs}

\title{Derived intersections over the Hochschild cochain complex}
\author{M\'arton Hablicsek\thanks{Mathematics Department, Leiden University,
Netherlands, {\em e-mail: }{\tt hablicsekhm@math.leidenuniv.nl}}}
\begin{document}
\maketitle
\input{abstract}
\input{intro}
\input{tdo}
\input{poly}
\input{dil}

\input{biblio}
\end{document}

%% file: defs.tex
\DeclareFontFamily{U}{rsf}{}
\DeclareFontShape{U}{rsf}{m}{n}{
  <5> <6> rsfs5 <7> <8> <9> rsfs7 <10-> rsfs10}{}
\DeclareMathAlphabet{\mathscr}{U}{rsf}{m}{n}

\DeclareMathAlphabet{\mathgth}{U}{euf}{m}{n}

\DeclareFontFamily{U}{cyr}{}
\DeclareFontShape{U}{cyr}{m}{n}{
  <5> wncyr5 <6> wncyr6 <7> wncyr7 <8> wncyr8 <9> wncyr9 <10-> wncyr10}{}
\DeclareMathAlphabet{\mathcyr}{U}{cyr}{m}{n}

\input cyracc.def

\makeatletter
\def\operator@font{\sf}
\makeatother

\newcommand{\sbt}{{\scalebox{0.5}{\textbullet}}}
\newcommand{\eps}{\epsilon}

\newcommand{\cA}{{\mathscr A}}
\newcommand{\cB}{{\mathscr B}}

\newcommand{\cC}{{\mathscr C}}

\newcommand{\cO}{{\mathscr O}}

\newcommand{\sTor}{\underline{\mathsf{Tor}}}

\newcommand{\sExt}{\underline{\mathsf{Ext}}}

\newcommand{\ra}{\rightarrow}

\newcommand{\C}{\mathbb{C}}

\newcommand{\iso}{\cong}

\newcommand{\del}{\partial}

\newcommand{\field}[1]{\mathbb{#1}}

\renewcommand{\L}{\field{L}}
\renewcommand{\phi}{\varphi}

%% file: abstract.tex
\begin{abstract}
Let $M$ and $N$ be Lagrangian submanifolds of a complex symplectic manifold $S$. A result of Behrend-Fantechi and Baranovsky-Ginzburg provides a Batalin-Vilkovisky structure on $\sExt^\sbt_{\cO_S}(\cO_X,\cO_Y)$, the dual of the structure sheaf of the derived intersection of the Lagrangian submanifolds. 

In this paper paper we generalize this result to the case of the 1-shifted cotangent bundle $T^*[1]X$ of a smooth scheme $X$ over $\C$. We show how one can obtain the twisted cotangent bundles as derived intersections of special Lagrangians in $T^*[1]X$, moreover we show that the derived intersection of the quantized Lagrangians coincide with the canonical quantization of the twisted cotangent bundles.
\end{abstract}

%% file: intro.tex
\section{Introduction}

\paragraph Let $S$ be a Poisson manifold over a field $k$ of characteristic 0 and let $P$ be a Poisson bivector on $S$. Consider a flat $k[\hbar]/(\hbar^2)$-deformation, $\cA$ of $\cO_S$ such that its commutator induces the Poisson bracket given by $P$.  Let $Y$ and $Z$ be two smooth coisotropic subvarieties of $S$ and assume that we have a flat $k[\hbar]/(\hbar)^2$-deformation, $\cB$ of $\cO_Y$ to a right $\cA$-module and a flat deformation, $\cC$ of $\cO_Z$ to a left $\cA$-module.  Associated to this data, greatly inspired by the work of Behrend and Fantechi (\cite{BehFan}), Baranovsky and Ginzburg (\cite{BarGin}) constructed a non-trivial second order differential operator 
\[\delta:\sTor_i^{\cO_S}(\cO_Y,\cO_Z)\ra \sTor_{i+1}^{\cO_S}(\cO_Y,\cO_Z)\]
(here $\sTor$ denotes the sheaf-Tor functor) which squares to 0.  This construction makes the algebra $\sTor_\sbt^{\cO_S}(\cO_Y,\cO_Z)$ a Batalin-Vilkovisky (BV for short) algebra.

\paragraph Let us explain the construction in the simplest case, in the case of the cotangent bundle of a smooth scheme $X$ over $\C$, $S=T^*X$ (where we consider the natural symplectic structure on $S$). In this case the canonical quantization of $\cO_S$ is convergent, at $\hbar=1$, we have the sheaf of differential operators, $D_X$. Assume further that $X$ is a Calabi-Yau variety over $k$ and let $f$ be a global regular function on $X$.  Consider the zero section, $Y=X$: $i:X\ra T^*X$ and the graph of $df$, $Z$ inside $T^*X$.  The Tor sheaves $\sTor_\sbt^{\cO_{T^*X}}(\cO_X,\cO_Z)$ are given by the cohomology of the derived tensor product
\[\cO_X\otimes^L_{\cO_{T^*X}} \cO_Z.\]
Given a volume form $\omega$ on $X$, we have a trivialization $\omega_X\iso \cO_X$. Using this trivialization $\omega_X\iso \cO_X$, the complex $\cO_X\otimes^L_{\cO_{T^*X}} \cO_Z$ is quasi-isomorphic (up to a shift) to the twisted de Rham complex 
\[0\ra \cO_X\xrightarrow{\wedge df}\Omega^1_X\xrightarrow{\wedge df} \Omega^2_X\ra...\]
Now, consider the ring of differential operators, $D_X$. The algebra $\cO_Z$ has a natural left $D_X$-module structure.  Moreover, given a volume form $\omega$ on $X$, since $\omega_X$ has a natural right $D_X$-module structure, we can regard $\cO_X$ as a right $D_X$-module as well. The derived tensor product
\[\omega_X\otimes^L_{D_X} \cO_Z\]
has a natural BV-structure.  Using the trivialization $\omega_X\iso \cO_X$ given by the volume form $\omega$ again, this derived tensor product (over $X$) is quasi-isomorphic (up to a shift) to twisted de Rham complex
\[0\ra \cO_X\xrightarrow{d+\wedge df} \Omega^1_X\xrightarrow {d+\wedge df} \Omega^2_X\ra...\]
In this case, the BV-differential is basically given by the de Rham differential.

\paragraph We can organize our spaces in the following diagram
\[\xymatrix{W\ar[r]\ar[d]& Y\ar[d]^0 \\ Z\ar[r]^{df} & T^*X}\]
where $W$ is the derived scheme given by the derived intersection of the zero section, $Y$ and the graph of $df$, $Z$ inside $T^*X$. The cotangent bundle has a canonical symplectic structure, the two embeddings, $Y\ra T^*X$ and $Z\ra T^*X$ have canonical Lagrangian structures (in the sense of \cite{PTVV}). Therefore, their derived intersection, $W$ is equipped with  a $-1$-shifted symplectic structure (\cite{PTVV}). In the Costello-Li framework (see, for instance, \cite{CosLi}), the quantizations of $\cO_W$ (where $W$ is equipped with a $-1$-shifted symplectic structure) is a BV-algebra structure on the sheaf $\cO_W[[\hbar]]$ which is, in this case, exactly the BV-algebra structure constructed above.

Summarizing the above discussion we see that the quantization of the derived intersection of Lagrangians inside the cotangent bundle is the derived intersection of the quantized Lagrangians. We can state a meta statement:\\

{\em The quantization of the derived intersection of Lagrangians is the derived intersection of the quantized Lagrangians.}
\medskip

\paragraph The purpose of the current paper is to give a further evidence of the meta statement by generalizing the above picture when the ambient space is the 1-shifted cotangent bundle $T^*[1]X$ (we assume that $X$ is a smooth scheme over $\C$.) We consider two special closed embeddings (in the sense of \cite{TV}): the zero section $Y=X\ra T^*[1]X$ and the graph of $\alpha$, $Z\ra T^*[1]X$, where $\alpha\in H^1(X,\Omega^1_X)$. These closed embeddings can be equipped with Lagrangian structures whenever $d\alpha\in H^1(X,\Omega^{\geq 2}_{X})=0$. (Here $\Omega^{\geq 2}_{X}$ denotes the truncation of the de Rham complex of $X$.)  We see that their derived intersection, $W$ is $0$-shifted (again, by \cite{PTVV}), moreover one can show that $W$ is a twisted cotangent bundle over $X$ (corresponding to $\alpha$). 
\[\xymatrix{W\ar[r]\ar[d]& X\ar[d]^0 \\ Z\ar[r]^{\alpha} & T^*[1]X}\]
(For more details see Section \ref{sec:tdo}.)
We remark that the symplectic structure on $W$ depends on the choices of the Lagrangian structures on $X\ra T^*[1]X$ and $Z\ra T^*[1]X$.

\paragraph In the framework of \cite{PTVV} the quantization of an $n$-shifted symplectic derived stack, $S$ is a deformation of the structure sheaf $\cO_S$ over $\C[[\hbar]]$ as a sheaf of $E_{n+1}$-algebras. In the case of $S=T^*[1]X$ with its natural 1-shifted symplectic structure, the canonical quantization of $\cO_{T^*[1]X}$ is convergent, at $\hbar=1$, we have the sheaf of the Hochschild cochain complex of $X$. The two closed embeddings $X\ra T^*[1]X$ and $Z\ra T^*[1]X$ are endowed with Lagrangian structures, the quantization of these Lagrangian morphisms are left and right module-structures on $\cO_X$ and $\cO_Z$ over the $E_2$-algebra of the Hochschild cochain compex of $X$. We provide such module structures in Proposition \ref{prop:quant}. Specifically, we consider the Hochschild cochain complex of $X$ as a brace algebra, and we equip the structure sheaf $\cO_X$ with both left- and right-brace module structures over the Hochschild cochain complex of $X$. In this way, $\cO_X$ has not just the structure of a non-commutative algebra, but it has the structure of left- or right-brace module structures over the Hochschild cochain complex as well.

\paragraph The key point is that once the left- and right-brace module structures are provided, the derived tensor product will have an \textbf{associative algebra structure}. Explicitly, we follow Safronov (\cite{Saf}) to equip the derived intersection of the quantizations of the Lagrangian morphisms $X\ra T^*[1]X$ and $Z\ra T^*[1]X$ with an $E_1$-algebra structure (\textbf{associative algebra structure}). First, we consider a specific global model of the Hochschild cochain complex which is equipped with a brace algebra structure (\cite{GerVor}). We equip the sheaves $\cO_X$ (and $\cO_Z$) with a natural structure of a left (and of a right) brace module over this brace algebra (see Proposition \ref{prop:quant}). Given a brace algebra, the derived tensor product of a left and a right brace module is endowed with an associative structure (\cite{GerVor}, \cite{Saf}). 

We are ready to state our main result (for a precise statement, see Theorem \ref{thm:main2}).

\begin{Theorem}\label{thm:main1}
The derived intersection (in the above sense) of $\cO_X$ and $\cO_Z$ over the Hochschild cochain complex is the canonical quantization of the associated twisted cotangent bundle, $W$.
\end{Theorem}
\medskip

\paragraph[Remark:] Given a symplectic 1-shifted derived Artin stack $S$ and two Lagrangian morphisms $X\ra S$ and $Y\ra S$, it is expected that the meta statement is true, i.e quantization of the derived intersection $W=X\times_S Y$ is given by the derived tensor product of the quantized Lagrangians (if exists) over the quantization of $S$ (see \cite{Pri} and \cite{Saf}). The current paper provies the first explicit evidence to this statement (or in general to meta statement  when the ambient derived Artin stack is an $n$-shifted symplectic derived stack for $n>0$).

\paragraph The paper is organized as follows. In Section \ref{sec:tdo} we review the notions of the rings of twisted differential operators and twisted cotangent bundles. In Section \ref{sec:poly} we explain how to regard the Hochschild cochain complex as a brace algebra (\cite{GerVor}), and then we endow the structure sheaves $\cO_X$ and $\cO_Z$ with left- and right-brace module structures (Propositiion \ref{prop:quant}). Finally, in Secton \ref{sec:dil} we present the proof of our main theorem, Theorem \ref{thm:main1}. First, we provide an associative algebra morphism from the canonical quantization of the twisted cotangent bundle to the derived tensor product of the brace modules $\cO_X$ and $\cO_Z$ over the Hochschild cochain complex of $X$; then we show that this morphism is a quasi-isomorphism.

\paragraph{\textbf{Acknowledgment.}} The author expresses his thanks to Jonathan Block, Damien Calaque, Tony Pantev, Jon Pridham, Pavel Safronov, Junwu Tu and Shilin Yu for useful conversations.

%% file: tdo.tex
\section{Twisted differential operators}
\label{sec:tdo}

In this section we review the notions of the rings of twisted differential operators and twisted cotangent bundles.  Our references are \cite{BeiKaz} and \cite{Gin}.  We also show how to obtain the twisted cotangent bundles as derived intersections of Lagrangians inside the shifted cotangent bundle $T^*[1]X$.

We begin with the definition of a ring of twisted differential operators.

\begin{Definition}
A sheaf of rings, $D$ of twisted differential operators (TDO for short) on $X$ is a positively filtered sheaf of associative algebras, whose associated graded, $gr D$ is isomorphic to the symmetric algebra of the tangent bundle:
\[gr D\iso Sym_{\cO_X} T_X.\]
\end{Definition}

\paragraph In particular, we have a short exact sequence
\[0\ra \cO_X\ra D^1\ra T_X\ra 0\]
where $D^1$ denotes the first filtered piece of the sheaf of twisted differential operators $D$. The TDO is completely determined by the so-called Atiyah-algebra structure on $D^1$. The set of TDO's can be classified as follows.

\begin{Theorem}\label{thm:tdo} (\cite{BeiKaz}, Lemma A.1.6) The set of TDO's is in one-to-one correspondence between elements of the cohomology group $H^1(X,\Omega^{\geq 1}_X)$.
\end{Theorem}
\medskip

(Here $\Omega^{\geq 1}_X$ is the complex  $\Omega^1_X\xrightarrow{d} \Omega^2_X\xrightarrow{d} \Omega^3_X...$).

\paragraph Let us briefly explain the above result. The exact triangle
\[\Omega^{\geq 2}_{X}[-1]\ra \Omega^{\geq 1}_X\ra \Omega^1_X\ra \Omega^{\geq 2}_{X}\]
induces a long exact sequence on cohomology:
\[...\ra H^0(X,\Omega^{\geq 2}_{X})\ra H^1(X,\Omega^{\geq 1}_X)\ra H^1(X,\Omega^1_X)\ra...\]
The image of the class of TDO inside $H^1(X,\Omega^1_X)$ is the class of the extension
\[0\ra \cO_X\ra D^1\ra T_X\ra 0\]
where we forget the Atiyah-algebra structure on the first filtered piece $D^1$ of $D$.
Moreover given an extension
\[0\ra \cO_X\ra D^1\ra T_X\ra 0\]
the Atiyah-algebra structures on $D^1$ form a torsor over the image of $H^0(X,\Omega^{\geq 2}_{X})$ inside $H^1(X,\Omega^{\geq 1}_X)$. 

We remark that the space $H^0(X,\Omega^{\geq 2}_{X})$ is the space of global closed 2-forms on $X$, $H^0(X, \Omega^2_{cl, X})$.

\paragraph Twisted differential operators are important objects in deformation quantization as they be obtained as canonical deformation quantizations of the structure sheaves of the twisted cotangent bundles. Let us explain this result briefly. We start with the definition of twisted cotangent bundles.

\begin{Definition}
A twisted cotangent bundle $W$ on $X$ is an $\Omega^1_X$-torsor $\pi:W\ra X$ such that
\begin{itemize}
\item The space $W$ is endowed with a symplectice form $\omega$,
\item The fibers of $\pi$ are Lagrangian with respect to this symplectic form,
\item For any 1-form $\nu$ on $X$, $t_\nu^*(\omega)=\pi^*(d\nu)+\omega$ where $t_\nu:W\ra W$ is the torsor-action corresponding to $\nu$. 
\end{itemize}
\end{Definition}
\medskip

Similarly to Theorem \ref{thm:tdo} one can classify the twisted cotangent bundles as well. We have

\begin{Theorem} (\cite{BeiKaz}, Lemma A.1.10)
The set of twisted cotangent bundles is in one-to-one correspondence between elements of the cohomology group $H^1(X,\Omega^{\geq 1}_X)$.
\end{Theorem}
\medskip

\paragraph\label{par:tdofontos} Let us explain the above result as well. Again, we consider the long exact sequence on cohomology
\[...\ra H^0(X,\Omega^{\geq 2}_{X})\ra H^1(X,\Omega^{\geq 1}_X)\ra H^1(X,\Omega^1_X)\ra...\]
induced by the exact triangle
\[\Omega^{\geq 2}_{X}[-1]\ra \Omega^{\geq 1}_X\ra \Omega^1_X\ra \Omega^{\geq 2}_{X}.\]
The image of the class of the twisted cotangent bundle in $H^1(X,\Omega^1_X)$ classifies $W$ (as a space) up to isomorphism as follows. Since $W\ra X$ is an $\Omega^1_X$-torsor, there is a covering $\{U_i\ra X\}$ such that the restriction of $W$ to each $U_i$ is isomorphic to the cotangent bundle $T^*U_i$. The gluing data gives rise to the class in $H^1(X,\Omega^1_X)$.

Once the class in $H^1(X,\Omega^1_X)$ is fixed, the set of twisted cotangent bundle structure on $W$ (the set of possible symplectic structures satisfying the required properties) forms a torsor over the image of $H^0(X,\Omega^{\geq 2}_{X})$ inside $H^1(X,\Omega^{\geq 1}_X)$.

\paragraph For a twisted cotangent bundle the pushforward of the structure sheaf, $\pi_*\cO_W$ carries a Poisson-algebra structure. The canonical quantization is convergent, at $\hbar=1$, we have the TDO corresponding to the class of the twisted cotangent bundle inside $H^1(X,\Omega^{\geq 1}_X)$. 

\paragraph Now, we show how to obtain the twisted cotangent bundles (with the symplectic structure) as the derived intersections of Lagrangian morphisms to the shifted cotangent bundle $T^*[1]X$.  Consider the two closed smooth subschemes of $T^*[1]X$ given by the zero section $X\ra T^*[1]X$ and the graph of $\alpha\in H^1(X,\Omega^1_X)$ which we denote by $Z$. We remark that as a scheme $Z$ is abstractly isomorphic to $X$, but they are not isomorphic as subschemes of $T^*[1]X$. 

The tangent complex of $T^*[1]X$ (over $X$) is quasi-isomorphic to $T_X\oplus \Omega^1_X[1]$, the cotangent complex $L_{T^*[1]X}$ (over $X$) is quasi-isomorphic to $T_X[-1]\oplus \Omega^1_X$. The underlying 2-form of the natural 1-symplectic structure on $T^*X[1]$ is given (at least over $X$) by the (shifted) identity morphism $id[1]:T_{T^*[1]X}\ra L_{T^*[1]X}[1]$ (see \cite{Cal} for more details).

We recall the definition of a Lagrangian morphism $Y\ra T^*[1]X$ in the sense of \cite{PTVV}. 

\begin{Definition} A smooth subscheme $i:Y\ra T^*[1]X$ is equipped with a Lagrangian structure if 
\begin{itemize} 
\item the pullback of the underlying 2-form $i^*id[1]\in \cA^2_{cl,Y}[1]$ is homotopic to 0, and
\item the homotopy induces a quasi-isoxmorphism between the tangent complex, $T_i$ of $i$ and the cotangent bundle $\Omega^1_Y$ of $Y$:
$T_i\ra \Omega^1_Y$. 
\end{itemize}
\end{Definition}

\noindent
We emphasize that a Lagrangian structure is a choice of a path between 1-shifted symplectic form and 0 in the space of closed closed two-forms $\cA^{2}_{cl}(X,1)$ inducing the above conditions.

\begin{Proposition}
The graph of $\alpha\in H^1(X,\Omega_X^1)$, $Z\ra T^*[1]X$ can be equipped with a Lagrangian structure if and only if $d\alpha\in H^1(X,\Omega^{\geq 2}_{X})$ is 0, where $d:\Omega^1_X\ra \Omega^{\geq 2}_{X}$ is the de-Rham differential.
\end{Proposition}
\medskip

\begin{Proof}
It is easy to see that the pullback form $i^*id[1]$ gives rise to the class of $d\alpha\in \pi_0(\cA^2_{cl}(X,1))=H^1(X,\Omega^{\geq 2}_{X})$ (using the isomorphism $Z\iso X$ and the fact that the 1-shifted symplectic form is induced by the Liouville 1-shifted 1-form). As a result, $d\alpha$ has to be 0 if the graph of $\alpha$ can be equipped with a Lagrangian structure. 

The other direction is proven in Theorem 2.22 in \cite{Cal}. More precisely, it is shown in \cite{Cal} that any lift $\bar{\alpha}\in H^1(X,\Omega^{\geq 1} _{X})$ of $\alpha$ provides a Lagrangian structure on the graph of $\alpha$, $Z$.\qed
\end{Proof}

\paragraph Consider the derived intersection $W$ of the graph of $\alpha\in H^1(X,\Omega^1_X)$ and the zero section:
\[\xymatrix{W\ar[r]\ar[d]& X\ar[d]^0 \\ Z\ar[r]^{\alpha} & T^*[1]X}.\]
Locally $\alpha$ is trivial, therefore the derived intersection of the zero section and the graph of $\alpha$ is locally isomorphic to the derived self-intersection of the zero section, which is the cotangent bundle, $T^*X$. The class $\alpha$ specifies the gluing data of these cotangent bundles, and we see that for a given $\alpha$, the derived intersection, $W\ra X$, is isomorphic to the twisted cotangent bundle (as a space) corresponding to $\alpha$. As a consequence, $\alpha$ determines the isomorphism class (as a space) of the twisted cotangent bundle.

On the other hand, the possible Lagrangian structures on $Z$ determine the symplectic structure on the twisted cotangent bundle. Indeed, the zero section $X\ra T^*[1]X$ has a canonical Lagrangian structure. Hence, given a Lagrangian structure on $Z\ra T^*[1]X$, the derived intersection of $X$ and $Z$ is endowed with a symplectic structure. 
The possible Lagrangian structures on the graph of $\alpha$ are parametrized by the first homotopy group $\pi_1(\cA^{2}_{cl})(X,1), \star)$ of the space of shifted closed 2-forms on $X$. This homotopy group is isomorphic to $H^0(X,\Omega^{\geq 2}_{X})$. As a consequence, we get an action of $H^0(X,\Omega^{\geq 2}_{X})$ on the possible symplectic structures on the twisted cotangent bundle as described in Paragraph \ref{par:tdofontos}.

%% file: poly.tex
\section{Polydifferential operators}\label{sec:poly}

In this section we define the complex $Diff(\cO_X^\sbt,\cO_X)$, which is a global model for the Hochschild cochain complex equipped with its brace algebra structure. Similarly, for every $\cO_X$-bimodule, $P$, we define the complex $Diff(\cO_X^\sbt,P)$, which is a global model of the Hochschild cochain complex with coefficients in $P$ equipped with its brace module structure (see \cite{CalvdB}, \cite{Yek} for more details). We also show how $Diff(\cO_X,D)$ is the convergent quantization of $\cO_Z$ at $\hbar=1$ (see \cite{Pri} for more details).

We begin with the definition of (poly)differential operators.

\begin{Definition}
Let $P$ be an $\cO_X$-bimodule, and $A:\cO_X\ra P$ a $k$-linear map. Given a sequence of functions $f_0, f_1, ..., \in \cO_X$, define a sequence of $k$-linear maps $A_m:\cO_X\ra P$ given by $A_{-1}=A$ and, $A_n:= f_n A_{n-1}-A_{n-1}f_n$. We say that $A$ is a differential operator of order at most $N$ if for every point $x\in X$ and every section $s\in \cO_X$ defined at $x$, there exists a neighborhood $U$ of $x$ and $N\geq 0$ such that for any open subset $V\subset U$ and any choice of functions $f_0, f_1,..., f_N$ on $V$, so that $A_N(s|_V)$ vanishes.
\end{Definition}
\medskip 

The differential operators $\cO_X \ra  P$ form a sheaf that we denote by $Diff_X(\cO_X, P)$. In the case of $P=\cO_X$ we denote the sheaf of differential operators by $D_X:=Diff_X(\cO_X,\cO_X)$.

\paragraph[Remark:] The sheaf of differential operators, $D_X$ is a TDO.

\begin{Definition} A $k$-polylinear map $A: \cO_X\times...\times \cO_X\ra P$ (of $n$ arguments) is a polydifferential
operator of (poly)order at most $(N_1,..., N_n)$ if it is a differential operator of order at most $N_j$ in the $j$-th argument whenever the remaining $n-1$ arguments are fixed.
\end{Definition}
\medskip

The polydifferential operators $\cO_X\times...\times \cO_X\ra P$ of $i$ arguments form a sheaf, which we denote by $Diff(\cO_X^i,P)$.

\paragraph We can identify the sheaf $Diff(\cO_X^i,P)$ with the tensor product
\[D_X\otimes_{\cO_X}....\otimes_{\cO_X}D_X\otimes_{\cO_X}P\]
where the number of the $D_X$ terms is $i$ as follows. The map
\[D_X\otimes_{\cO_X}....\otimes_{\cO_X}D_X\otimes_{\cO_X}P\ra Diff(\cO_X^i,P)\]
given by
\[A_1\otimes...\otimes A_i\otimes p\mapsto A_1(-)A_2(-)...A_i(-)p\]
(for local sections $A_i\in D_X$ and $p\in P$) is clearly an isomorphism. (Here we use the natural $\cO_X$-bimodule structure on $D_X$.) We denote the isomorphism $Diff(\cO_X,P)\iso D_X\otimes P$ by $i_P$.

\paragraph The sheaves of polydifferential operators form a natural complex $Diff(\cO_X^\sbt,P)$ whose $i$-th term is $Diff(\cO_X^i,P)$ and the differential $d$ is given by 
\begin{align*}
dA(g_1,...,g_{i+1})&=g_1A(g_2,...,g_{i+1})-A(g_1g_2,g_3,...,g_{i+1})+A(g_1,g_2g_3,...,g_{i+1})+\\
&+(-1)^iA(g_1,...,g_ig_{i+1})+(-1)^{i+1}A(g_1,...,g_i)g_{i+1}
\end{align*}
where $A:\cO_X\times...\times \cO_X\ra P$ is a polydifferential operator of $i$ arguments and $\{g_1,...,g_{i+1}\}$ is a local section of $\cO_X\times...\times \cO_X$ (of $i+1$ arguments).

\paragraph Let $A\in Diff(\cO_X^i,\cO_X)$ and $B\in Diff(\cO_X^j,\cO_X)$. We define $A\cdot B\in Diff(\cO_X^{i+j},\cO_X)$ as the differential operator mapping $a_1,...,a_{i+j}$ to 
\[(-1)^{ij}A(a_1,...,a_i)\cdot B(a_{i+1},...,a_{i+j}).\]
This product endows the complex $Diff(\cO_X^\sbt,\cO_X)$ a differential graded algebra structure.

\paragraph Similarly, we see that if $P$ is an $\cO_X$-bimodule, then the complex $Diff(\cO_X^\sbt,P)$ has a differential graded bimodule structure over $Diff(\cO_X^\sbt,\cO_X)$. Moreover, if the $\cO_X$-bimodule, $P$ has an associative algebra structure, then $Diff(\cO_X^\sbt,P)$ has a differential graded algebra structure defined parallel to the differential graded algebra structure on $Diff(\cO_X^\sbt,\cO_X)$.

\paragraph The complex, $Diff(\cO_X^\sbt,\cO_X)$ is equipped with a brace algebra structure as follows. Let $A\in Diff(\cO_X^i,\cO_X)$ and $A_{l}\in Diff(\cO_X^{j_l},\cO_X)$ (for $l=1,...,m$). The brace operations $A\{A_1,...,A_m\}$ are defined as operations of degree $-m$, i.e. $A\{A_1,...,A_m\}\in Diff(\cO_X^n,\cO_X)$, where $n=i+\sum_{l=1}^m j_l-m$. Explicitely, $A\{A_1,...,A_m\}$ maps $n$ sections of $\cO_X$, $a_1,...,a_n$ to
\[\sum_{0\leq i_1\leq...\leq i_m\leq n}(-1)^\epsilon A(a_1,...,a_{i_1},A_1(a_{i_1+1},...),...,a_{i_m},A_m(a_{i_m+1},...),...,a_n)\]
where $\eps:=\sum_{i=1}^m i_l(j_l-1)$. 

Similarly, if $P$ is an $\cO_X$-bimodule, then $Diff(\cO_X^\sbt,P)$ has a module structure over the brace algebra $Diff(\cO_X^\sbt,\cO_X)$ as follows. Let $B\in Diff(\cO_X^\sbt,P)$ and $A_{l}\in Diff(\cO_X^{j_l},\cO_X)$ (for $l=1,...,m$). The brace operations $B\{A_1,...,A_m\}$ are defined as operations of degree $-m$, i.e. $B\{A_1,...,A_m\}\in Diff(\cO_X^n,P)$, where $n=i+\sum_{l=1}^m j_l-m$. Explicitly, $B\{A_1,...,A_m\}$ maps $n$ sections of $\cO_X$, $a_1,...,a_n$ to
\[\sum_{0\leq i_1\leq...\leq i_m\leq n}(-1)^\epsilon B(a_1,...,a_{i_1},A_1(a_{i_1+1},...),...,a_{i_m},A_m(a_{i_m+1},...),...,a_n)\]
where $\eps:=\sum_{i=1}^m i_l(j_l-1)$. 

\paragraph A right $Diff(\cO_X^\sbt,\cO_X)$ brace module structure on a complex is equivalent to a left $Diff(\cO_X^\sbt,\cO_X)^{op}$ brace module structure. We equip the complex $Diff(\cO_X^\sbt,D_X^{op})^{op}$ with a right $Diff(\cO_X^\sbt,\cO_X)$ structure, coming from the left $Diff(\cO_X^\sbt,\cO_X)$ brace module structure on $Diff(\cO_X^\sbt,D_X^{op})$.

\paragraph Let $P=D$ be a TDO. Then the complex $Diff(\cO_X^\sbt,D)$ starts as 
\[D\xrightarrow{d_D} Diff(\cO_X,D)\ra...\]
where the differential acts as
\[d_D(\del)=(x\mapsto \del x-x\del).\]
Clearly, the kernel of this map is the centralizer of $\cO_X$ inside $D$, which is $\cO_X$. Hence, the zeroth cohomology sheaf of the complex $Diff(\cO_X^\sbt,D_X)$ is $\cO_X$. Moreover, all other cohomology sheaves vanish (see \cite{Sri} for instance), hence we obtain the following quasi-isomorphism of algebra objects.

\begin{Proposition}\label{prop:quant} We have
\[Diff(\cO_X^\sbt,D)\cong \cO_X\]
(and similarly
\[Diff(\cO_X^\sbt,D^{op})^{op}\cong \cO_X).\]
\end{Proposition}

%% file: dil.tex
\section{Derived intersection of the quantizied Lagrangians}
\label{sec:dil}

In this section we prove our main theorem, Theorem \ref{thm:main1}. We begin with a simple but key lemma.

\begin{Lemma}\label{le:dil}
Let $D$ be a TDO. The image of the map
\[\phi:D\ra Diff(\cO_X,D)\iso D_X\otimes_{\cO_X} D\]
given by 
\[\del\ra i_D(d_D\del)+1\otimes \del\]
is the kernel of the map
\[Diff(\cO_X, D)\xrightarrow{d_{Diff(\cO_X,D)}} Diff(\cO_X, Diff(\cO_X, D)).\]
\end{Lemma}

\begin{Proof}
First, the map $\phi: D\ra Diff(\cO_X,D)$ can be realized as the map
\[\del \mapsto (x\mapsto \del \cdot x)\]
where $\del\cdot x $ is given by the multiplication on $D$. In other words, $\phi(\del)$ is just the map given by multiplication from the left by $\del$. 

Moreover, an element $f\in Diff(\cO_X, D)$ is in the kernel of $d_{Diff(\cO_X,D)}$ if and only if for every local section $x\in \cO_X(U)$ we have that $fx=xf$. This means that for any local section $y\in \cO_X(U)$ we have that $f(y)x=f(yx)$ meaning that $f$ has to be given by left multiplication by an element $\del \in D$.\qed
\end{Proof}
\medskip

\paragraph As a consequence of the Lemma above, we obtain a morphism 
 \[D\ra Diff(\cO_X^{\sbt}, Diff(\cO_X, D)).\]
The cup product map 
\[\cup: Diff(\cO_X^\sbt, D_X^{op})^{op}\otimes_{Diff(\cO_X^\sbt, \cO_X)}Diff(\cO_X^\sbt, D)\ra Diff(\cO_X^\sbt, D_X\otimes^\L_{\cO_X} D)\]
provides the quasi-isomorphism between the right-hand side and the derived tensor product
\[Diff(\cO_X^\sbt, D_X^{op})^{op}\otimes_{Diff(\cO_X^\sbt, \cO_X)}Diff(\cO_X^\sbt, D).\]
Composing with this quasi-isomorphism, we obtain a morphism of complexes
\[\chi: D\ra Diff(\cO_X^\sbt, D_X^{op})^{op}\otimes_{Diff(\cO_X^\sbt, \cO_X)}Diff(\cO_X^\sbt, D).\]

\paragraph \label{par:rem} Before we prove our main theorem  we recall how the associative algebra structure is defined on
\[Diff(\cO_X^\sbt, D_X^{op})^{op}\otimes_{Diff(\cO_X^\sbt, \cO_X)}Diff(\cO_X^\sbt, D).\]
We follow (\cite{Saf}). Let $A$ be a brace algebra and $N$ (and $M$ resp.) be a right (and left resp.) brace module over $A$. The bar complex $T_\sbt(A[1])$ is a dg coalgebra. Gerstenhaber and Voronov \cite{GerVor} defined a multiplicative structure on the bar complex which makes $T_\sbt(A[1])$ into a bialgebra. Similarly, the one-sided bar complexes $N\otimes T_\sbt(A[1])$ and $T_\sbt(A[1])\otimes M$ carry a natural structure of a dg algebra compatible with the $T_\sbt(A[1])$ comodule structure. Therefore, the cotensor product 
\[N\otimes T\sbt(A[1]) \otimes^{T_\sbt(A[1])}T_\sbt(A[1])\otimes M\]
carries a natural structure of a dg algebra. This cotensor product is quasi-isomorphic to the derived tensor product $N\otimes_A^L M$, providing the required dg algebra structure on the derived tensor product. This dg algebra structure is compatible with the dg-algebra structure on $N\otimes_\C M$ meaning that the natural map 
\[N\otimes_\C M\ra N\otimes_A^\L M\]
is a dg-algebra map (\cite{Saf}).

\paragraph Recall that in Section \ref{sec:poly} we equipped $Diff(\cO_X^\sbt, \cO_X)$ with a brace algebra structure, and moreover, we provided left- and right-brace module structures on $Diff(\cO_X^\sbt, D)$ and $Diff(\cO_X^\sbt, D_X^{op})^{op}$. As a result, we obtain a dg-algebra structure on 
\[Diff(\cO_X^\sbt, D_X^{op})^{op}\otimes_{Diff(\cO_X^\sbt, \cO_X)}Diff(\cO_X^\sbt, D).\]

We are ready to prove our main theorem.

\begin{Theorem}\label{thm:main2}
The derived tensor product
\[Diff(\cO_X^\sbt,D_X^{op})^{op}\otimes_{Diff(\cO_X^\sbt,\cO_X)} Diff(\cO_X^\sbt,D)\]
is quasi-isomorphic to $D$ as associative algebras (under the quasi-isomorphism $\chi$).
\end{Theorem}

\begin{Proof}
By Lemma \ref{le:dil} and the discussion followed we obtain a map of complexes 
\[\chi: D\ra Diff(\cO_X^\sbt,D_X^{op})^{op}\otimes_{Diff(\cO_X^\sbt,\cO_X)} Diff(\cO_X^\sbt,D).\]
First, we show that this map is an algebra map. We will prove this by showing that \'etale locally the algebra structure on $D$ can be lifted to $Diff(\cO_X^\sbt,D_X^{op})^{op}\otimes_\C Diff(\cO_X^\sbt,D)$. Since, we are in a local setting, we assume that $D=D_X$ and $X$ is an affine $n$-space over $\C$. In other words, we assume that $D=D_X=\C\langle x_1,...,x_n, \del_1,...,\del_n\rangle$ where the relations are given by the usual Weyl-algebra relations.

Conider the following elements in $D_X\otimes_\C D_X$:
\[\bar{\del}_i:=\del_i\otimes 1+1\otimes \del_i, \quad\mbox{and}\quad \bar{x}_j:=1\otimes x_j\]
where both $i$ and $j$ run from 1 to $n$. These elements, satisfy the usual Weyl-algebra relations
\[\bar{\del}_i\bar{x}_j-\bar{x}_j\bar{\del}_i=\delta_{ij}.\]
Using these elements, there is a unique algebra map $\rho: D_X\rightarrow D_X\otimes_\C D_X$ sending the $\del_i$ to the $\bar{x}_i$ and the $x_j$ to the $\bar{x}_j$. Under the projection 
\[\pi: D_X\otimes_\C D_X\ra D_X\otimes_{\cO_X} D_X\]
these elements become the corresponding elements under the map $\phi: D_X\ra Diff(\cO_X, D_X)\cong D_X\otimes_{\cO_X} D_X$. Indeed, a general element $\del \in D_X$ given as $\del=f\prod_{i\in I} \del_i$, we have
\[(\pi\circ\rho)(\del)=\sum_{J\subset I}\prod_{j\in J}\del_j \otimes f\prod_{i\in I\setminus J}\del_i.\]
Similarly,
\[\phi(\del)=(x\mapsto f\prod_{i\in I} \del_i \cdot x)=\sum_{J\subset I} (\prod_{j\in J} \del_j).x\cdot \prod_{i\in I\setminus J}\del_i\] showing that $(\pi\circ \rho)(\del)=\phi(\del)$.

Since  $D_X\iso Diff(\cO_X^0,D_X^{op})^{op}$ and $D\iso Diff(\cO_X^0, D)$ as algebras, we obtain the corresponding elements $\bar{\del}_i$ and $\bar{x}_j$ inside $Diff(\cO_X^0,D_X^{op})^{op}$ and $D\iso Diff(\cO_X^0, D)$. By the discussion at the end of Paragraph \ref{par:rem} we obtain that the map
\[\chi: D\ra Diff(\cO_X^\sbt,D_X^{op})^{op}\otimes_{Diff(\cO_X^\sbt,\cO_X)} Diff(\cO_X^\sbt,D)\]
is a dg-algebra map!

Finally, we show that $\chi$ is a quasi-isomorphism. First, we note that $Diff(\cO_X^\sbt, D_X^{op})^{op}$ has a natural filtration given by the filtration of $D_X$. This filtration is compatible with the $Diff(\cO_X^\sbt, \cO_X)$-action. This induces a filtration on the derived tensor product 
\[Diff(\cO_X^\sbt,D_X^{op})^{op}\otimes_{Diff(\cO_X^\sbt,\cO_X)} Diff(\cO_X^\sbt,D).\]
However, by Proposition \ref{prop:quant}, we know that $Diff(\cO_X^\sbt,D)$ is quasi-isomorphic to $\cO_X$, in other words, the graded pieces of the fitration are all concentrated in degree 0. As a consequence, the cohomology of the complex 
\[Diff(\cO_X^\sbt,D_X^{op})^{op}\otimes_{Diff(\cO_X^\sbt,\cO_X)} Diff(\cO_X^\sbt,D)\]
has to be concentrated in degree 0 as well. Hence, by Lemma \ref{le:dil}, we are done.\qed
\end{Proof}

\paragraph[Remark:] The proof above implies that the map of complexes
\[\chi: D\ra Diff(\cO_X^\sbt,D_X^{op})^{op}\otimes_{Diff(\cO_X^\sbt,\cO_X)} Diff(\cO_X^\sbt,D)\]
is a quasi-isomorphism of \textit{filtered associative algebras}. 

\medskip

%% file: project1.bbl
\begin{thebibliography}{GP}

\bibitem{AriBloPan} Arinkin, D., Block, J., Pantev, T., {\em *-Quantizations of Fourier-Mukai Transforms}, Geometric and Functional Analysis, (2013), Volume 23, Issue 5, pp 1403--1482

\bibitem{BarGin} Baranovsky, V., Ginzburg, V., {\em Gerstenhaber-Batalin-Vilkoviski structures on coisotropic
intersections}, Mathematical Research Letters 17.2 (2010): 211--229.

\bibitem{BehFan} Behrend, K., Fantechi, B., {\em Gerstenhaber and Batalin-Vilkovisky structures on Lagrangian
intersections}, in Algebra, Arithmetic and Geometry: Manin Festschrift, Volume I., Vol. 269
of Progress in Mathematics, Birkh\"auser (2009)

\bibitem{BeiKaz} Beilinson, A., Kazhdan, D., {\em Flat Projective Connections}, unpublished

\bibitem{Cal} Calaque, D., {\em Shifted cotangent stacks are shifted symplectic}, Annales de la Faculté des sciences de Toulouse: Mathématiques. Vol. 28. No. 1. (2019)

\bibitem{CalvdB} Calaque, D., Van den Bergh, M., {\em Hochschild cohomology and Atiyah classes}, Adv. Math. 224
(2010), no. 5, 1839--1889.

\bibitem{CosLi} Costello, K., Li, S., {\em Quantum BCOV theory on Calabi-Yau manifolds and the higher genus $\beta$-model}, preprint, (2012)


\bibitem{GerVor} Gerstenhaber, M., Voronov, A., {\em Higher operations on the Hochschild complex}, Func. Anal. and its
App. 29 (1995) 1--5.

\bibitem{Gin} Ginzburg, V., {\em Lectures on $D$-modules}, lecture notes


\bibitem{PTVV} Pantev, T., To\"en, B., Vaqui\'e, M., Vezzosi, G., {\em Shifted symplectic structues}, Publ. math. de l'IHÉS,  (2013), Volume 117, Issue 1, pp 271--328 

\bibitem{Pri} Pridham, J., P.; {\em Quantization of derived lagrangians}, preprint, (2016)

\bibitem{Saf} Safronov, P., {\em Poisson reduction as a coisotropic intersection}, Higher Structures 1.1 (2018).

\bibitem{Sri} Sridharan, R., {\em Filtered Algebras and Representations of Lie Algebras}, Trans. Am. Math. Soc., Vol. 100, No. 3 (1961), pp. 530--550

\bibitem{TV} To\"en, B., Vezzosi G., {\em Homotopical algebraic geometry. II. Geometric stacks and applications},
Memoirs of the American Mathematical Society, American Mathematical Society, (2008), 193
(902)

\bibitem{Yek} Yekutieli, A., {\em The continuous Hochschild cochain complex of a scheme}, Canad. J. Math. 54
(2002), no. 6, 1319--1337.

































\end{thebibliography}
